\documentclass[10pt,a4paper]{article}
\usepackage[utf8]{inputenc} 
\usepackage{amsmath}
\usepackage{amsfonts}
\usepackage{amssymb}
\newcommand{\zzzP}{\mathbb{Z}'}
\newcommand{\Mo}{\mathbb{M}}
\usepackage{epsfig}
\usepackage{graphics}
\DeclareMathAlphabet\mathbfcal{OMS}{cmsy}{b}{n}
\newcommand{\zbT}{\mathbfcal{T}}
\newcommand{\ZA}{\mathcal{A}}
\newcommand{\ZLL}{\langle\!\langle}
\newcommand{\ZRR}{\rangle\!\rangle}
\newcommand{\ZSUno}{\sum _{n=1}^{+\ZIN}}

\newcommand{\zln}{\lambda_n}
 
\newcommand{\zb}[1]{{\bf #1}}
\newcommand{\zbg}[1]{{\boldsymbol #1}}

\newcommand{\zL}{\mathcal{L}}

\newcommand{\zbn}{\zb{n}}
 
  \newcommand{\zbPsi}{\zbg{\Psi}}
 
  \newcommand{\zbv}{\zb{v}}
 \newcommand{\zbu}{\zb{u}}
  \newcommand{\zbf}{\zb{f}}
   \newcommand{\zbgg}{\zb{g}}
 \newcommand{\zbphi}{\zbg{\phi}}
\newcommand{\zbw}{\zb{w}}
\newcommand{\zbpsi}{\zbg{\psi}}

\newcommand{\zbxi}{\zbg{\xi}}
\newcommand{\zbeta}{\zbg{\eta}}
\newcommand{\zbF}{\zb{F}}

\newcommand{\intt}{\int_0^t}
\newcommand{\intT}{\int_0^T}
\newcommand{\ints}{\int_0^s}
\newcommand{\intr}{\int_0^r}

\newcommand{\ZOMq}{\Omega}

\newcommand{\zl}{\lambda}

\newcommand{\ZDE}{\delta} 

\newcommand{\zdia}{~~\rule{1mm}{2mm}\par\medskip}  

\newcommand{\ZIN}{\infty}
\newcommand{\zProof}{{\bf\underbar{Proof}.}\ } 
  
\newcommand{\zg}{\gamma} 
\newcommand{\zaa}{\alpha} 
\newcommand{\zzr}{{\rm I\hskip-2.1pt R}}

\newcommand{\ZBI}{\bibitem}

\newcommand{\ZD}{\;\mbox{\rm d}}    
  
\newcommand{\ZEP}{\epsilon}

\newcommand{\ZLA}{\label}

\newcommand{\ZSI}{\sigma}

\newtheorem{Theorem}{Theorem}  
\newtheorem{Corollary}[Theorem]{Corollary}  
\newtheorem{Lemma}[Theorem]{Lemma} 

\newtheorem{Remark}[Theorem]{Remark}
\newtheorem{Definition}[Theorem]{Definition} 
\newtheorem{Example}[Theorem]{Example}

\author{
L. Pandolfi\thanks{Dipartimento di Scienze Matematiche ``Giuseppe Luigi Lagrange'', Politecnico di Torino, Corso Duca degli Abruzzi 24, 10129 Torino, Italy (luciano.pandolfi@polito.it)}
}
\title{Controllability of  isotropic viscoelastic bodies of Maxwell-Boltzmann type\thanks{
This papers fits into the research program of the GNAMPA-INDAM and has been written in the framework of the   ``Groupement de Recherche en Contr\^ole des EDP entre la France et l'Italie (CONEDP-CNRS)''.}}

\begin{document}

\maketitle

 {\bf\underline{Abstract}:} In this paper we consider a viscoelastic   three dimensional body (of Maxwell-Boltzmann type)  controlled on (part of) the boundary. We assume that the material is isotropic and homogeneous.  If the body is elastic (i.e. no dissipation due to past memory), controllability has been studied by several authors. We prove that the viscoelastic body inherits the controllability properties of the corresponding purely elastic system. The proof relays on cosine operator methods combined with moment theory.

 \section{Introduction}
 
 We consider a linear viscoelastic body occupying a region $\ZOMq\subseteq \zzr^3$. We assume that the body is isotropic
 and homogeneous so that the dynamics of the body is described by the operator $ \zL $:
 \[ 
 \zL\zbu=\mu\Delta \zbu +(\zl+\mu)\nabla\left ( \nabla\cdot\zb u\right )\,,\quad \zbu=\zbu(x)\quad x\in\ZOMq\,.
  \]
 Here $\zbu\in\zzr^3 $ and  $ \zl >0$,  $ \mu>0 $ are the   Lam\'e coefficients. We assume that $\zl$ and $\mu$ are constant (i.e. that the body is homogeneous).

Boldface denotes vectors. Note that also the space variable (i.e. $ x $) is a vector, with $ {\rm dim}\, x={\rm dim}\, \zbu $, but we don't use boldface for the space variable (and not for the   vector $ 0 $).
 
If the body is elastic, the evolution in time of the displacement is
described by the Navier  equation
\begin{equation}\ZLA{equa:puroELASICO}
\zbu''=\zL \zbu+\zbF\,,\qquad \left\{\begin{array}{l}
\zbu(0)=\zbu_0 \\
\zbu'(0)=\zbu_1 
\end{array}\right.
\end{equation}
($\zbu=\zbu(x,t)$, $  \zbF=\zbF(x,t) $ and $x\in\ZOMq\subseteq \zzr^3$) while if the body is viscoelastic (of the Maxwell-Boltzmann type) the evolution in time of the displacement (denoted $\zbw=\zbw(x,t)\in\zzr^3$)  is described by  
 \begin{equation}\ZLA{equa:viscoELASICO}
\zbw''=\zL\zbw +\intt M(t-s)\zL\zbw(s) \ZD s+\zbF\,,\qquad \left\{\begin{array}{l}
\zbw(0)=\zbw_0 \\
\zbw'(0)=\zbw_1  
\end{array}\right.
\end{equation}
(see~\cite[p.~112-113]{Kolsky}) In fact, the integral should extend from $-\ZIN$ but in the study of reachability for linear systems, the control being applied after a time $t_0$, say $t_0=0$,   it is not restrictive to assume $\zbw=0$ for $t<t_0=0$.  
 
 We   apply a displacement on a part $\Gamma$ of  the boundary of $\ZOMq$ ($\Gamma=\partial\ZOMq$ is not excluded), both to the elastic and the viscoelastic body, and we are going to study  whether is it possible to control the pair ((\emph{velocity of displacement})/(\emph{displacement}))  to hit a prescribed target  at a certain time $T$. 
 I.e. we impose the following boundary condition for $t>0$: 
 \begin{eqnarray} 
 \ZLA{condiBOrdoELASTICO}
  \zbu(x,t)  =\zbf(x,t) \  x\in \Gamma  & &  \zbu(x,t)  =0 \  x\in \partial\ZOMq\setminus\Gamma \\
   \ZLA{condiBOrdoViscoELASTICO}  
\zbw(x,t) = \zbf(x,t)\  x\in \Gamma  & &\!\! \ \zbw(x,t) = 0 \ x\in \partial\ZOMq\setminus\Gamma  
  \end{eqnarray}
  (if $\Gamma=\partial\ZOMq$ then disregard the condition on $\partial\ZOMq\setminus\Gamma$).
   
 When the initial condition and the affine term $ \zbF $ are zero and we want to stress dependence on the control $\zbf$, we write  $\zbu^\zbf$, $\zbw^\zbf$. 
  
Under the assumptions we shall state below, it turns out that $ (\zbu'(t),\zbu(t)) $ and $ (\zbw'(t),\zbw(t)) $ are $\left [ H^{-1}(\ZOMq)\right ]^3\times  \left [L^2(\ZOMq)\right ]^3 $-valued continuous functions so that their values at a fixed time $T$ make sense in this space. 
 
  If  the initial conditions are zero and also $ \zbF=0 $, the \emph{reachable sets  at time $T$ } are the sets  \begin{align*}
  R_E(T)&= \left\{(\left (\zbu^\zbf\right )'(T),\zbu^\zbf(T))\,, \ \; \zb f\in L^2\left (0,T;\left [L^2(\Gamma)\right ]^3\right )\right \} \,,\\
R_V(T)&= \left\{(\left (\zbw^\zbf\right )'(T),\zbw^\zbf(T))\,, \  \zb f\in L^2\left (0,T;\left [L^2(\Gamma)\right ]^3\right ) \right \} \,. 
\end{align*}   
  Note that (both in the purely elastic   and in the viscoelastic case) the reachable set increases with time.

\emph{Controllability at time $ T $} is the property   $  R_{V}(T)=   \left [ H^{-1}(\ZOMq)\right ]^3\times \left [L^2(\ZOMq)\right ]^3 $ for the viscoelastic system, and 
$  R_{E}(T)=  \left [ H^{-1}(\ZOMq)\right ]^3\times \left [L^2(\ZOMq)\right ]^3 $ in the elastic case.

 Under the assumptions we describe below, controllability in the purely elastic case, i.e.  for the system~(\ref{equa:puroELASICO})-(\ref{condiBOrdoELASTICO}), has been studied  by several authors (see the references in 
 Section~\ref{sec:References}). We are going to prove that the controllability property which holds in the purely elastic case is inherited by the viscoelastic system.
 
\subsection{Notations, assumptions and the main results of this paper}
In this section we state the assumptions, we describe preliminary results on the controllability of the elastic system and we state our main results, which will be proved in the next sections.

{\bf\underline{Assumption A)}}
We assume:
\begin{itemize}
 \item
the kernel $M(t)$ is of class $H^2_{\rm loc}([0,+\ZIN))$.  
 
\item   the region $\ZOMq\subseteq \zzr^3$  is bounded and   $ \partial\ZOMq $ is of class $C^2$.

 \item We assume that both the Lam\'e constants $ \zl $ and $ \mu  $ are positive and constant (i.e. the body is homogeneous). 
 
  \item The subset $\Gamma$ of $\partial\ZOMq$ will be called the \emph{active part} of the boundary. The first assumption on $\Gamma $ is that it is relatively open in $\partial\ZOMq$.

  \end{itemize}

Known facts for the elastic system~(\ref{equa:puroELASICO}) with boundary
 conditions~(\ref{condiBOrdoELASTICO}), prved in the references cited in Sect.~\ref{sec:References}: 
\begin{itemize}
\item Let $T>0$ and $\zb f=0$. Let $\zbF\in L^1\left (0,T;\left [L^2(\ZOMq)\right ]^3\right )  $,     $\zbu_0\in [H^1_0(\ZOMq)]^3$, $\zbu_1\in [L^2(\ZOMq)]^3$. Then 
\[
t\mapsto (\zbu(t),\zbu'(t))\in C\left ([0,T];[H^1_0(\ZOMq)]^3\right )\times C\left ([0,T];[L^2(\ZOMq)]^3\right )\,.
\]

\item Let $T>0$. Let $\zbF\in L^1\left (0,T;[L^2(\ZOMq)]^3\right )$, $f\in L^2\left (0,T ;[L^2(\Gamma)]^3\right )$ and $\zbu_0\in [L^2(\ZOMq)]^3$, $\zbu_1\in [H^{-1}(\ZOMq)]^3$. Then 
\[
t\mapsto (\zbu(t),\zbu'(t))\in C\left ([0,T];[L^2(\ZOMq)]^3\right )\times C\left ([0,T];[H^{-1}(\ZOMq)]^3\right )\,.
\]

\item  
\emph{There exist open subsets $\Gamma$  of $\partial\ZOMq$ and times $T$ such that    the purely elastic system   is controllable in time $T$, i.e. 
\begin{equation}\ZLA{controEq:Ela}
 R_{E}(T)=  \left [ H^{-1}(\ZOMq)\right ]^3\times \left [L^2(\ZOMq)\right ]^3\,.
\end{equation}
 }
\end{itemize}

\noindent{\bf\underline{Assumption B)}} The active part $\Gamma$ of $\partial\ZOMq$ and $T >0$ are chosen in such a way that the purely elastic system~(\ref{equa:puroELASICO}) with boundary control~(\ref{condiBOrdoELASTICO})  is controllable in time $T$, i.e. equality~(\ref{controEq:Ela}) holds.

 \begin{Remark}[On the notations]{\rm 
 From  now on, for the sake of readability, we   drop the exponent $3$ and we    write simply $L^2(\ZOMq)$, $H^1_0(\ZOMq)$,  $H^{-1}(\ZOMq)$ instead of $\left [L^2(\ZOMq)\right ]^3$,   $\left [H^{1}_0(\ZOMq)\right ]^3$, $\left [H^{-1}(\ZOMq)\right ]^3$.}
  \end{Remark}

We recall the following integration by parts formula, which we shall repeatedly use:
  
 \begin{equation}\ZLA{eq:INTEperPARTI}
 \int _{\ZOMq}\left ( \mathcal{L}\zbu\right )  \cdot\zbphi\ZD x=
\int_\Gamma  \left (\mathcal{T}\zbu\right )\cdot\zbphi\ZD\Gamma 
 -\int_\Gamma \left (\mathcal{T}\zbphi\right ) \cdot\zbu\ZD\Gamma+\int_\ZOMq \zbu\cdot \left (\mathcal{L}\zbphi\right )\ZD x 
 \,.
 \end{equation}
 The boundary operator $\zbT$, the \emph{boundary traction,}   is:
\[
\zbT\zbphi=
\mu\left (\zbn\cdot\nabla\right )\zbphi +(\zl+\mu)\left (\nabla\cdot \zbu\right )\zbn
\]
where  $\zbn$ is the exterior normal to $\partial\ZOMq$.

The following property, known as (output) \emph{admissibility} or \emph{direct inequality} is proved for example in~\cite[Ch.~IV]{Lions}:
\begin{Lemma}
Let $\zbu$ solve~(\ref{equa:puroELASICO}) and~(\ref{condiBOrdoELASTICO}) with $\zbf=0$. For every $T>0$ there exists a number $M $  such that
\begin{equation}
\ZLA{DIREelastico}
\begin{array}{l}
|\mathcal{T}\zbu|_{L^2( 0,T;L^2(\Gamma))}\leq M\left (  |\zbu_0|^2_{H^1_0(\ZOMq)}+|\zbu_1|^2_{L^2(\ZOMq)} +|\zbF|^2_{L^1(0,T;L^2(\ZOMq))}\right )\,.
\end{array}
\end{equation}
\end{Lemma}

Note that the inequality usually proved is  
\[
 |\zg_1\zbu|_{L^2( 0,T;L^2(\Gamma))}\leq M\left (  |\zbu_0|^2_{H^1_0(\ZOMq)}+|\zbu_1|^2_{L^2(\ZOMq)}+|\zbF|^2_{L^1(0,T;L^2(\ZOMq))} \right ) 
\] 
($\zg_1$ is the normal derivative on $\partial\ZOMq$) from which~(\ref{DIREelastico}) follows because $\zbu=0$ on $\partial\ZOMq$ implies $\nabla  u_k=\zg_1 u_k$ for every component $u_k$ of $\zbu$.

 Our first, ancillary result is as follows.  
 \begin{Theorem}\ZLA{teo:risuGene}
Let Assumption A) hold. Then:
 \begin{enumerate} 
\item\ZLA{item1inteo:risuGene} Let $ \zbF\in  L^1(0,T;L^2(\ZOMq))  $ and $\zbf=0$.
 System~(\ref{equa:viscoELASICO}) with initial conditions $\zbw(0)=\zbw_0\in  H^1_0(\ZOMq) $, $\zbw'(0)=\zbw_1\in L^2(\ZOMq) $ 
 admits a unique solution $\zbw\in C\left ([0,T];  H^1_0(\ZOMq) \right )\cap 
 C^1\left ([0,T]; L^2(\ZOMq) \right )$ for every $T>0$.
 \item\ZLA{item2teo:risuGene} Let $ \zbF\in  L^1(0,T;L^2(\ZOMq))  $.
 System~(\ref{equa:viscoELASICO}) with initial conditions $\zbw(0)=\zbw_0\in  L^2(\ZOMq) $, $\zbw'(0)=\zbw_1\in  H^{-1}(\ZOMq) $ and boundary   control $ \zb f\in  L^2(0,T;L^2(\Gamma)) $
 admits a unique solution $\zbw\in C\left ([0,T];  L^2(\ZOMq) \right )\cap 
 C^1\left ([0,T]; H^{-1}(\ZOMq) \right )$ for every $T>0$.
 \item\ZLA{item3teo:risuGene}   Let $\zbf=0$, $\zbw_0\in H^1_0(\ZOMq)$,  $\zbw_1\in L^2(\ZOMq)$ and let $\zbw$ 
 solve~(\ref{equa:viscoELASICO})  and~(\ref{condiBOrdoELASTICO}). Then, for every $T>0$ and every $\Gamma\subseteq\partial\ZOMq$ there exists $M $ such that
 \begin{equation}\ZLA{DireINEQ}
|\mathcal{T}\zbw|_{L^2( 0,T;L^2(\Gamma))}\leq M\left (  |\zbw_0|^2_{H^1_0(\ZOMq)}+|\zbw_1|^2_{L^2(\ZOMq)}
+|\zbF|^2_{L^1(0,T;L^2(\ZOMq))} \right )\,.
 \end{equation}
 This inequality is  the \emph{direct inequalities} of Eq.~(\ref{equa:viscoELASICO}).
 \end{enumerate}
 \end{Theorem}
 
 We noted already that the statement of Theorem~\ref{teo:risuGene} holds when $M=0$.
 
 The control result which we intend to prove is:
 \begin{Theorem}\ZLA{teo:controllabilita}
  Let Assumptions A)  and  B) hold so that
$
 R_{ E}(T)=  H^{-1}(\ZOMq) \times L^2(\ZOMq) 
$
and let $T_1>T$.
Then, system~(\ref{equa:viscoELASICO}),~(\ref{condiBOrdoViscoELASTICO})  is controllable in time $T_1$, i.e.  $R_{ V}(T_1)= H^{-1}(\ZOMq) \times  L^2(\ZOMq)  $.
 
 \end{Theorem}

 The proof is based on two main steps. {\bf Step 1)} proves that if the purely elastic system~(\ref{equa:puroELASICO})  is controllable at time $T$   then the reachable set $R_{ V}(T)$ of the viscoelastic system  is closed with finite codimension in $H^{-1}(\ZOMq)\times L^2(\ZOMq)$. Thanks to the fact that the reachable set is increasing, this holds at every $T_1\geq T$.    {\bf Step 2)} proves  that if the purely elastic system is controllable at time $T $ and $T_1>T$  then we have
 \[
 R_{ V}(T_1)= H^{-1}(\ZOMq )\times L^2(\ZOMq )\,.
 \]

 The organization of the paper is as follows: the next subsection recalls  previous results while section~\ref{sect:Cosine} presents  preliminaries on the cosine operator theory and the proof of Theorem~\ref{teo:risuGene}.

The two steps of the proof of Theorem~\ref{teo:controllabilita} are in Section~\ref{teo:controllabilita}.

 \subsection{\ZLA{sec:References}References to previous work and preliminaries}
 
Controllability of the Navier equation~(\ref{equa:puroELASICO}) has been studied in several papers, of which we cite~\cite{YamaGRASSE,YamaGRASSE2,Lions,Telega}. The proofs in these papers  are based on the \emph{inverse inequality,} obtained using multiplier methods, i.e. it is proved that if $\Gamma$ and $T $ are suitably chosen then there exists $m>0$ such that
\begin{equation}\ZLA{inVEinERAuNO}
m\left (  |\zbw_0|^2_{H^1_0(\ZOMq)}+|\zbw_1|^2_{L^2(\ZOMq)}
 \right )\leq  \intT\int_\Gamma\left [\mu |\zg_1\zbw|^2+(\zl+\mu)|\nabla\cdot\zbw|^2\right]^2\ZD\Gamma\,\ZD t\,. 
\end{equation}
Here it is assumed $\zbf=0$, $\zbF=0$ and of course $\zbw_0\in H^1_0(\ZOMq)$, $\zbw_1\in L^2(\ZOMq)$.
As noted in~\cite[p.~228]{Lions}, this inequality implies 
\begin{equation}\ZLA{inVEinERAdUE}
m\left (  |\zbw_0|^2_{H^1_0(\ZOMq)}+|\zbw_1|^2_{L^2(\ZOMq)}
 \right )\leq 
\intT\int_\Gamma \left |\mathcal{T}\zbw\right |^2\ZD\Gamma\,\ZD t 
\end{equation}
(for a different $m>0$) and this is the ``dual''  formulation  of the definition of  controllability.

 We mention that the inverse inequality  in the paper~\cite{YamaGRASSE2}   is proved also if the Lam\'e coefficients  are (slowly) space varying and that isotropy is not assumed    in~\cite{Telega}. These results can be extended to the viscoelastic case, precisely as we do below in the isotropic homogeneous case.

 In this paper, systems~(\ref{equa:puroELASICO}) and~(\ref{equa:viscoELASICO}) are studied using a cosine operator approach and also moment methods. 
  Let us define the following operators    in $L^2(\ZOMq)$:
  \begin{equation}
\ZLA{eq:ch1:defiOperA}
{\rm dom}\, A=H^2(\ZOMq)\cap H^1_0(\ZOMq)\;,\qquad A\zbphi=\mathcal{L}\zbphi\,, \quad \ZA=i (-A)^{1/2} \;.
\end{equation}
 The operator $A$ is selfadjoint with compact resolvent. Hence, $L^2(\ZOMq)$ admits an orthonormal basis of eigenvectors of $A$. Let it be denoted $\{\zbphi_n\}$ and let $-\zl_n^2$ be the eigenvalue of $\zbphi_n$ (it is known that the eigenvalues are negative and we assume that they are ordered in such a way that $\zl_n\leq \zl_{n+1}$).
 The operator $\ZA$ generates a $C_0$-group of operators on $\zzr$ so that we can define the operators
 \[
 R_+(t) =\frac{1}{2}\left[
e^{\ZA t}+e^{-\ZA t}
\right]\,,\qquad R_-(t) =\frac{1}{2}\left[
e^{\ZA t}-e^{-\ZA t}
\right]
 \]
 (the operator $R_+(t)$ is called the \emph{cosine operator generated by $A$).}  
 \emph{Note that    $t\mapsto R_+(t)\zbu$ and $t\mapsto R_-(t)\zbu$ belong to $\mathcal{L}\left (L^2(\ZOMq),C([0,T],L^2(\ZOMq))\right )$ and also to $\mathcal{L}(H^1_0(\ZOMq),C\left ([0,T],H^1_0(\ZOMq))\right )$.}
 
The operators $R_+(t)$ and $R_-(t)$ have the following expansions in series of $\{\zbphi_n\}$
\begin{equation}
\ZLA{ExpansCosINESINE}
\begin{array}{l}
\displaystyle R_+(t)\left (
\ZSUno\zaa_n\zbphi_n
\right ) =\ZSUno \left (\zaa_n\cos \zl_n t\right ) \zbphi_n(x)\;,
\\
\displaystyle
 R_-(t)\left (
\ZSUno\zaa_n\zbphi_n
\right )
 =i\left (\ZSUno \left (\zaa_n\sin \zl_n t\right ) \zbphi_n(x)\right ) \;.
\end{array} 
\end{equation}
 The solutions of problem~(\ref{equa:puroELASICO})-(\ref{condiBOrdoELASTICO}) are given by (see~\cite{BieliLASIE})
 
 \begin{equation}\ZLA{eq:SoLUsoLoElasT}
\zbu(t)=R_+(t)\zbu_0+\ZA^{-1} R_-(t)\zbu_1+\ZA^{-1}\intt R_-(t-s)\zbF(s)\ZD s-\ZA\intt R_-(t-s)D\zb f(s)\ZD s\,.
 \end{equation} 
 The operator $D$ in this formula is the operator $\zbf\mapsto \zbu=D\zbf$ where $\zbu$ solves
 \[
 \mathcal{L} \zbu=0\,,\quad  \left\{\begin{array}{l}
  \zbu(x,t)  =\zbf(x,t) \  x\in \Gamma\,,\\  \zbu(x,t)  =0 \  x\in \partial\ZOMq\setminus\Gamma\,.
 \end{array}\right.
 \]
 
It is known (see~\cite[Th.~3.6]{Dahlbert})  that $D\in\mathcal{L}\left (L^2(\partial\ZOMq),L^2(\ZOMq)\right )$ and, as noted in~\cite[p.~796]{Dahlbert},  it takes values in $H^{1/2}(\ZOMq)$, hence it is a compact operator.  

We use repeatedly the fact that when $\zbw\in  {\rm dom}\,A$ we have
 \begin{equation}\ZLA{eq:FormuFRONTiera}
 \zbT\zbw=-D^*A\zbw
 \end{equation} (a fact that can be proved using the  integration by parts~(\ref{eq:INTEperPARTI}) as in~\cite[p.~181]{LasieTRIvol1}).
 
 The cosine operator approach to controllability of systems with persistent memory, when $\zbu=u\in\zzr$, $\zbw=w\in\zzr$ and   $ \mathcal{L}=\Delta$,
  was first used in~\cite{PandAMO} where the \emph{existence} of the control time was proved, but the control time itself was not identified. The control time was identified in dimension
  $d=1$   in subsequent papers and using moment methods, see for example~\cite{AB,LoretiPANDOLFIsforza,PandIEOT,Pandcina} and for the scalar valued wave equation in $\ZOMq
  \subseteq \zzr^3$
  in~\cite{PandolfiSHARP}. In this paper we combine the cosine operator method and the moment method in order 
  to get a  proof of controllability for Eq.~(\ref{equa:viscoELASICO}).  This is based on the use of the following known estimates for the eigenvalues of the operator $A$ (see~\cite{ChengSTIMASI,Pleijel}): there exist $m>0$ and $M$ such that
 \begin{equation}\ZLA{eq:STImasintAUTOV}
m\left (n^{2/3}\right )\leq \zl_n^2\leq M\left (n^{2/3}\right )\,.
 \end{equation}

 \section{\ZLA{sect:Cosine}Cosine operators and   the proof of Theorem~\ref{teo:risuGene}}
  
 Now we prove Theorem~\ref{teo:risuGene} in the case $M\neq 0$.
The proof of the items~\ref{item1inteo:risuGene} and~\ref{item2teo:risuGene} is similar to that   in~\cite[Chapt.~2]{PandLIBRO}   and it is only sketched. 
  The first step is a definition of the solutions of problem~(\ref{equa:viscoELASICO})-(\ref{condiBOrdoViscoELASTICO}).  
 We first apply formally the MacCamy trick. Let $R(t)$ be the resolvent kernel of $M(t)$, given by  
\begin{equation}\ZLA{eq:resolvM}
R(t)+\intt M(t-s)R(s)\ZD s=M(t)\;.
\end{equation} 
 We ``solve'' the Volterra integral equation~(\ref{equa:viscoELASICO})
 in the ``unknown'' $\zL\zbw $. 
We get
\[
\zL\zbw(t)=\zbw''(t)-\zbF(t)-\intt R(t-s)[ \zbw''(s)-\zbF(s)]\ZD s\,.
\]

 We integrate by parts and we get   
 
\begin{equation}\ZLA{eq:foPoMacCam}
\left\{\begin{array}{l}
\displaystyle   \zbw'' = \zL \zbw+a\zbw'(t)  +b \zbw(t)+\intt K(t-s) \zbw(s)\ZD s +\zbF_1(t)\,,\\
\displaystyle  \zbF_1(t) = \zbF(t)+\intt R(t-s)\zbF(s)\ZD s-R(t)\zbw_1-R'(t)\zbw_0 \,,\\[2mm]
\displaystyle   a= R(0)\,, \quad 
 b=R'(0)\,,\quad K(t)=R''(t) 
\end{array}\right. 
 \end{equation}
 and the initial and boundary conditions~(\ref{condiBOrdoViscoELASTICO}).
 i.e. MacCamy trick removes the differential operator from the memory term.
 
 \begin{Remark}\ZLA{rema:peraequalzero}{\rm
 The term $a\zbw'(t)$ can be removed from the right hand side if we perform the transformation $\zbv(t)=e^{-(a/2)t}\zbw(t)$. 
The effect on the initial conditions is   that $\zbv(0)=\zbv_0=\zbw_0$ while $\zbv'(0)=\zbv_1=\zbw_1-(a/2)\zbw_0$ and the control   is replaced by $e^{-(a/2)t}\zbf(t)$. Of course this has no influence on controllability and so with obvious  modifications of the definitions  of $\zbf(t)$, $\zbF_1(t)$,  $b$ and $K(t)$ we shall work with~(\ref{eq:foPoMacCam}) but with $a=0$ (we shall use the expression $\zbv'(0)=\zbw_1-(a/2)\zbw_0$ when needed for clarity).  
 Moreover note    that the kernel $K(t)$ is square integrable.\zdia}
 \end{Remark}

We explicitly note that if  $(\zbw_0,\zbw_1)\in H^1_0(\ZOMq)\times L^2(\ZOMq)$ (respectively      $(\zbw_0,\zbw_1)\in   L^2(\ZOMq)\times H^{-1}(\ZOMq)$)  then
 $\zbF_1(t)$ belongs to $L^1(0,T;L^2(\ZOMq))$ (respectively to $L^1(0,T;H^{-1}(\ZOMq))$)   and it depends continuously on $(\zbw_0,\zbw_1)$ in the specified spaces.
  
  Using formula~(\ref{eq:SoLUsoLoElasT}) we get the following Volterra integral equation of the second kind for $\zbw(t)$:
  \begin{equation} \ZLA{eq:Cap1:VoltePERlasoluzioneWpostBIS}
  \zbw(t)  =  \zbu(t) + \ZA^{-1}\intt R_-( t-s) \zbF_1(s)\ZD s+ 
   L*\zb w 
  \end{equation}
 where $*$ denotes the convolution, 
  \begin{equation}\ZLA{eq:DefiDiL}
L(t)\zbw=   bR_-(t)\zb w+\intt K(t-r)R_-(r)\zbw\ZD r 
 \end{equation}
and 
 $\zbu(t)$ is the solution of~(\ref{equa:puroELASICO}) when $\zbu_0=\zbw_0$ and $\zbu_1=\zbw_1-(a/2)\zbw_0  $ and $\zbF(t)=0$ (recall   $a=0$, see Remark~\ref{rema:peraequalzero}):
 \[
\zbu(t)=- \ZA\intt R_-( t-s )D\zb f(s)\ZD s + R_+( t)\zbu_0+ \ZA^{-1}R_-( t)\zbu_1\,. 
 \]  
 Now we observe that the function
 \[
t\mapsto  \ZA^{-1}\intt R_-( t-s) \zbF_1(s)\ZD s 
 \]
  belongs to $C([0,T];H^1_0(\ZOMq))\cap C^1([0,T];L^2(\ZOMq))$  if $(\zbw_0,\zbw_1)\in H^1_0(\ZOMq)\times L^2(\ZOMq)$ and to
 $C([0,T];L^2(\ZOMq))\cap C^1([0,T];H^{-1}(\ZOMq))$  if $(\zbw_0,\zbw_1)\in  L^2(\ZOMq)\times H^{-1}(\ZOMq)$ and this is precisely the regularity of $\zbu(t)$. The fact that the solution of a Volterra integral equation 
 retains the regularity of the affine term proves items~\ref{item1inteo:risuGene} and~\ref{item2teo:risuGene} in Theorem~\ref{teo:risuGene}. Now we prove the direct inequality, i.e. the statement~\ref{item3teo:risuGene}.
 In this proof,  $T>0$ is fixed and we use $M$ to denote constants (possibly depending on $T$) which are  not the same at every occurrence.
 
 \begin{Remark}\ZLA{rema:perORToGOn}
 {\rm
 We shall use the following observation: the initial conditions enters into the affine term (and $\zbw_1$ is transformed to $\zbw_1-(a/2)\zbw_0$) when we used MacCamy trick. If the original equation to be studied is
\[
\zbpsi''=\zL \zbpsi +b\zbpsi(t)+\intt K(t-s)  \zbpsi(s) \ZD s \,,\qquad   \left\{\begin{array}{l}
 \zbpsi(0)= \zbxi \\
 \zbpsi'(0)= \zbeta \\  \zbpsi=0\ {\rm on}\ \partial\ZOMq\,.
 \end{array}\right.
\]
 then $\zbpsi$ solves the Volterra integral equation~(\ref{eq:Cap1:VoltePERlasoluzioneWpostBIS}) with $\zbF_1=0$ and $\zbu(t)=R_+(t)\zbxi+\ZA^{-1}R_-(t)\eta$.\zdia
 }
 \end{Remark}
 
 In the proof of the statement~\ref{item3teo:risuGene} of Theorem~\ref{teo:risuGene} we assume $\zbf=0$, $\zbw_0\in H^1_0(\ZOMq)$, $\zbw_1\in L^2(\ZOMq)$ and $\zbF\in L^1(0,T;L^2(\ZOMq))$.
 Let
 \begin{align*}
\tilde\zbu(t)&=\zbu(t) +\ZA^{-1} \intt R_-( t-s) \zbF_1(s)\ZD s=\\
&=
R_+(t)\zbu_0+\ZA^{-1}R_-(t)\zbu_1+\ZA^{-1} \intt R_-( t-s) \zbF_1(s)\ZD s
\,,\\
& \tilde\zbu(0)=\zbu_0=\zbw_0\in H^1_0(\ZOMq)\,,\quad \tilde \zbu'(0)=\zbu_1=\zbw_1-\frac{a}{2} \zbw_0\in L^2(\ZOMq)\,.
 \end{align*}
 Picard iteration applied to~(\ref{eq:Cap1:VoltePERlasoluzioneWpostBIS})  gives 
 \begin{equation}\ZLA{eq:Picard}
\zbw(t)=\tilde\zbu(t)  +\ZA^{-1}\intt L(t-s)\tilde \zbu(s)\ZD s   +A^{-1}\left [\sum _{n=2}^{+\ZIN} \ZA^{-n+2}\left (L^{(*n)}\right )*\tilde\zbu\right ](t) 
 \end{equation} 
 where $L^{(*n)}  $ denotes iterated convolution.

 We prove the direct  inequality in the case  $\zbw_0\in \mathcal{D}(\ZOMq) $, $\zbw_1\in  \mathcal{D}(\ZOMq)$ and  $\zbF\in \mathcal{D}(\ZOMq\times(0,T))$. 
In this case, $\zbw\in C([0,T];H^1_0(\ZOMq)\cap H^2(\ZOMq))=C([0,T];{\rm dom}\, A)$.  
 The inequality is then extended to $\zbw_0\in H^1_0(\ZOMq) $,  $\zbw_1\in L^2(\ZOMq)$,  $\zbF\in L^1(0,T;L^2(\ZOMq))$ by continuity and density.

Using~(\ref{eq:FormuFRONTiera}) we see that there exists $M=M(T)$ such that
 
 \begin{align*}
|  \zbT \zbw |^2_{L^2((0,T);L^2(\Gamma))}
  \leq &  M 
  \biggl (
  |\zbT \tilde \zbu |^2_{L^2((0,T);L^2(\Gamma))}   +   \biggr.\\
 & \biggl. 
  +\left | 
 \zbT\left (\ZA^{-1}\intt L(t-s)\tilde \zbu(s)\ZD s \right )
 \right|^2_{L^2(0,T;L^2(\Gamma))}   + \biggr.\biggr. \\
 &
   \biggl.+\left |D^*\left (\sum _{n=2}^{+\ZIN} \ZA^{-n+2}\left (L^{(*n)}\right )*\tilde\zbu \right )
 \right |^2_{L^2(0,T;L^2(\Gamma))} 
  \biggr )\,. 
 \end{align*}
 The known properties of the solution $\zbu$ of the elastic system~(\ref{equa:puroELASICO}) implies that (for a possibly different constant $M=M_T$)
 
 \begin{align*}
 |\zbT \tilde \zbu |^2_{L^2((0,T);L^2(\Gamma))}+\left |D^*\left (\sum _{n=2}^{+\ZIN} \ZA^{-n+2}\left (L^{(*n)}\right )*\tilde\zbu \right )
 \right |^2_{L^2(0,T;L^2(\Gamma))}\leq  \\   \leq M\left (|\zbw_0|^2_{H^1_0(\ZOMq)}+|\zbw_1|^2_{L^2(\ZOMq)}+|\zbF|^2_{L^1(0,T;L^2(\ZOMq))}\right)\,.
 \end{align*}
In order to complete the proof we show that we have also
\begin{align*}
&\left | 
 \zbT\left (\ZA^{-1}\intt L(t-s)\tilde \zbu(s)\ZD s \right )
 \right|^2_{L^2(0,T;L^2(\Gamma))}\leq\\  & \leq
 M\left (|\zbw_0|^2_{H^1_0(\ZOMq)}+|\zbw_1|^2_{L^2(\ZOMq)}+|\zbF|^2_{L^1(0,T;L^2(\ZOMq))}\right)\,.
\end{align*} 
 We note that
 \[
 \ZA^{-1}\intt L(t-s)\tilde u(s)\ZD s=
\ZA^{-1}\intt L(t-s) u(s) \ZD s+ A^{-1} \intt L(t-s)\ints R_-(s-r)\zbF_1(r)\ZD r\,\ZD s 
 \]
 and the required inequality holds for the second addendum.
 The first addendum is
\begin{align*}
&\ZA^{-1}\intt L(t-s)  u(s)\ZD s=\\
  &=\ZA^{-1}\left [\intt bR_-(t-s) R_+(s)\zbu_0\ZD s+\intt K(r)\int_0^{t-r} R_-(t-r-s) R_+(s)\zbu_0\ZD r\,\ZD s\right ]+\\
  &+ A^{-1}\left [
  \intt bR_-(t-s)  R_-(s)\zbu_1 \ZD s+\intt K(r)\int_0^{t-r} R_-(t-r-s)   R_-(s)\zbu_1\ZD r\,\ZD s
  \right ]\,.
\end{align*}
 
Using~(\ref{eq:FormuFRONTiera}), we see that:
 \begin{align*}
& \left |\mathcal{T} A^{-1}\left [
 \intt bR_-(t-s)   R_-(s)\zbu_1 \ZD s+\intt K(r)\int_0^{t-r} R_-(t-r-s)   R_-(s)\zbu_1\ZD r\,\ZD s
 \right ]\right |\leq \\
  & M |\zbu_1|_{L^2(\ZOMq)}\leq 
M\left (|\zbw_0|_{L^2(\ZOMq)}+|\zbw_1|_{L^2(\ZOMq)}\right ) 
\leq  
  M\left (   |\zbw_0|_{H^{1}_0(\ZOMq)}+|\zbw_1|_{L^2(\ZOMq)}    \right )
 \end{align*}
 (recall $\zbu_1=\zbw_1-(a/2)\zbw_0$).
 
We consider the term containing $\zbw_0$ using
\[
R_-(s)R_+(r)=\frac{1}{2} \left (R_-(s+r)+R_-(s-r)\right )\,.
\] 
So we have
\begin{align*}
\mathcal{T} \ZA^{-1}\intt R_-(t-s)R_+(s)\zbw_0\ZD s= \frac{t }{2}\mathcal{T}\ZA^{-1}R_-(t)\zbw_0+\frac{1}{2}D^*R_-(t)\zbw_0\,.
\end{align*}
 The second addendum, as an element of $C([0,T];L^2(\Gamma))$,  depends continuously on $\zbw_0\in H^1_0(\ZOMq) $.
 
 The function $\ZA^{-1}R_-(t)\zbw_0$ is the solution of~(\ref{equa:puroELASICO}) with $\zbf=0$, $\zbF=0$ and initial data $\zbu_0=0$, $\zbu_1=\zbw_0\in H^1_0(\ZOMq)\subseteq L^2(\ZOMq)$ so that
 \[
| \mathcal{T}\ZA^{-1}R_-(t)\zbw_0|^2_{L^2(0,T;L^2(\Gamma))}\leq M|\zbw_0|^2_{L^2(\ZOMq)} \leq M|\zbw_0|^2_{H^1_0(\ZOMq)}\,.
 \]
 The remaining iterated integral is treated analogously.

 See~\cite{PandMADRID} for similar arguments and~\cite{PandAMO} for a different proof (in the case of the wave equation).

\section{The proof of controllability}

In this section we prove the controllability result in Theorem~\ref{teo:controllabilita}. As noted, the proof is in two steps. Let $R_E(T)=H^{-1}(\ZOMq)\times L^2(\ZOMq)$. Then we prove
\begin{itemize}
\item   $R_V(T) $, at the same time $T$ (and so also at larger times), is closed with finite   codimension.
\item if $\ZEP>0$ and $\zbphi\perp R_V(T+\ZEP) $ then $\zbphi=0$ and so $R_V(T+\ZEP)=H^{-1}(\ZOMq)\times L^2(\ZOMq)$
\end{itemize}

\subsection{The first step}

The first step is simple and it is based on the Picard formula~(\ref{eq:Picard}). Let $\zbw_0=0$, $\zbw_1=0$ and $\zbF=0$ so that $\tilde \zbu=\zbu$.   We consider the operator $\zbf\mapsto\zbw^\zbf (T)$ which, according to~(\ref{eq:Picard}),   is the sum of three terms, which are linear and continuous as functions of $\zbf\in L^2(0,T;L^2(\Gamma))$. The first one
is  
\[
\zbf\mapsto \zbu(T)=  \zbu^{\zbf}(T)\,:
\]
this operator is surjective by assumptions (in fact, even more: $ \zbf  \mapsto \left (  \left (\zbu^\zb f\right )'(T),\zbu^\zbf(T) \right ) $ is surjective). The third operator 
\[
\zbf\mapsto  A^{-1}\left [\sum _{n=2}^{+\ZIN} \ZA^{-n+2}\left (L^{(*n)}\right )*\zbu\right ](T) 
\]
takes values in ${\rm dom}\, A=H^1_0(\ZOMq)\cap H^2(\ZOMq) $, and so it is a compact operator from   $  L^2(0,T;L^2(\Gamma))$ to $ L^2(\ZOMq) $ .

The operator $D$ is compact and so the operators of $\zbf\in L^2(0,T;L^2(\Gamma))$
\[
\zbf\mapsto \ZA^{-1} \intT L(T-s)\zbu^\zbf(s)\ZD s
=-\intT L(T-s)\ints R_-(s-r) D\zbf(r)\ZD r\,\ZD s
\]
with values in $L^2(\ZOMq)$  is compact.

The function $\zbf \mapsto \zbw'(T)$, with values in $H^{-1}(\ZOMq)$  is treated analogously. 
Note that 
\[
\frac{\ZD}{\ZD t}\ZA^{-1}L(t)\zbu=bR_+(t)\zbu+\intt K(r)R_+(t-r)\zbu\ZD r
\]
belongs to $C([0,T];L^2(\ZOMq))$ and $L^2(\ZOMq)$ is compactly embedded in $H^{-1}(\ZOMq)$ and so the transformation 
\[
\zbf \mapsto\left [\frac{\ZD}{\ZD t}\ZA^{-1}\intt L(t-s)\zbu(s)\ZD s\right ]_{t=T}
\] from $L^2(0,T;L^2(\Gamma))$ to $H^{-1}(\ZOMq)$ is compact.
The same holds for 
the derivative of the last term in~(\ref{eq:Picard}) computed for $t=T$, which is

\begin{align*}
&\ZA^{-1}\intt R_+(t-r)\left [ b\zbgg(r)+\intr K(t-r-s) \zbgg(s)\ZD s\right ]\ZD r\,,\\
&\zbgg =\zbu+ \sum _{n=3}^{+\ZIN } \ZA^{-n+2}L^{(*(n-1))}*\zbu\,.
\end{align*}

 In conclusion,
the reachable set $R_V(T)$ is the image of   the operator $\zbf\mapsto \left (\left ( \zbu^\zbf\right ) '(T), \zbu^\zbf(T)\right )$, which is surjective in $H^{-1}(\ZOMq)\times L^2(\ZOMq)$ since we assumed that    system~(\ref{equa:puroELASICO})  is controllable at time $T$, perturbed by the addition of a compact operator. This implies that  $R_V(T)$ is closed with finite codimension.  

\subsection{The second step}
The second step is more involved and requires several substeps.  

\subsubsection{Step~2-1: the elements of $\left  [R_V(T) \right ]^\perp$}
 
 We elaborate on formula~(\ref{eq:Picard}) which can be written as
 
\begin{equation}
 \ZLA{eq:PicardCOMPATTA}
 \zbw(t)=\tilde  \zbu(t)+\intt H(s)\tilde \zbu(t-s)\ZD s \,,\qquad H(t)\zbv=\ZSUno \ZA^{-n}L^{(*n)}  \zbv\,.
 \end{equation}

We study the reachable set, hence we assume $\zbF=0$ and null  initial conditions, so that

\[
\tilde\zbu(t)=\zbu(t)=\ZA\zbv(t)\,,\qquad 
\zbv(t) =\intt R_-(t-s)D\left (-\zbf(s)\right )\ZD s  \,.
\]
Of course, in the study of the reachable set, $-\zb f$ can be renamed $\zb f$.

We characterize the elements $(\zbxi,\zbeta)\in H^1_0(\ZOMq)\times L^2(\ZOMq) $ which annihilates $R_V(T)$ (here $T>0$ is arbitrary). 
I.e. we characterize those elements $(\zbxi,\zbeta)\in H^1_0(\ZOMq)\times L^2(\ZOMq)$ such that
\[
\ZLL\left (\zbw^\zbf\right )'(T),\zbxi\ZRR+\int_\ZOMq \zbw^\zbf(x,T)\zbeta(x)\ZD x=0
\]
for every $\zbf\in L^2(0,T;L^2(\Gamma))$. The crochet denotes the pairing of $H^{-1}(\ZOMq)$ and $H^1_0(\ZOMq)$. This set of annihilators is shortly denoted $\left  [R_V(T)\right ]^\perp$ and  we write $(\zbxi,\zbeta)\perp R_V(T)$.

Note that if it happens that  $\zbw'(T)\in L^2(\ZOMq)$ then 
\[
\ZLL\zbw '(T),\zbxi\ZRR=\int_\ZOMq \zbw'(x,T)\zbxi(x)\ZD x\,.
\]

Of course we can study  $  R_V(T)^\perp $ by assuming $f\in \mathcal{D}(\Gamma\times(0,T))$ so that the following computations are justified.  First we note that when $f\in \mathcal{D}(\Gamma\times(0,T))$ we have
$\zbw(t)\in L^2(\ZOMq)$ since
\[
\zbw'(t)=\ZA\intt R_-( s)D\zbf'(t-s)\ZD s+\ZA\intt H( s)\int_0^{t-s} R_-( r)D\zbf'(t-s-r)\ZD r\,\ZD s \,.
\] 
Hence (in the last step we integrate by parts in time and we note that $R_-(t)$ and $H(t)$ commute),
\begin{align}
\nonumber&\ZLL\zb w'(T),\zbxi\ZRR 
 =\int _{\ZOMq}\zbxi(x)\ZA\intT R_-(T-s)Df'(s)\ZD s\,\ZD x+\\
\nonumber &+\int_{\ZOMq}\zbxi(x)\ZA\intT H(s)\int_0^{T-s}R_-(T-s-r)Df'(r)\ZD r\,\ZD s\,\ZD x=\\
\nonumber&=\int_\Gamma\intT f'(s)\left [D^*\ZA R_-(T-s)\zbxi\right ]\ZD s\,\ZD\Gamma+\\
\nonumber&+\int_\Gamma\intT f'(r)\left [D^*\ZA \int_0^{T-r} H(s)R_-(T-r-s) \zbxi\ZD s\right ]\ZD r\,\ZD\Gamma=\\
\ZLA{SecCOMPO}&=\int_\Gamma\intT f(r)\left \{D^*A\left [R_+(T-r)\zbxi+\int_0^{T-r} H(s)R_+(T-r-s) \zbxi\ZD s\right ]\right \}\ZD r\,\ZD\Gamma\,.
\end{align}
 Note that the last integral makes sense thanks to the direct inequality,
   because~(\ref{eq:PicardCOMPATTA}) shows that the bracket  is the  solution of~(\ref{equa:viscoELASICO}) with $\zbf=0$, $\zbF=0$ and initial condition $\zbw(0)=\zbxi\in H^1_0(\ZOMq)$, $\zbw'(0)=0$.
 
 Analogously,
 \begin{align}
\nonumber &\int_\ZOMq \zbeta(x) \left [
 \ZA \intT R_-(T-r) D \zbf(r)\ZD r+\ZA\intT H( s)\int_0^{T-s} R_-( r)D \zbf(T-s-r)\ZD r\,\ZD s
\right ]\ZD x=\\
\nonumber&= \intT\int_\Gamma \zbf(r) D^*\ZA R_-(T-r)\zbeta\ZD r\,\ZD\Gamma+\intT\int_\Gamma \zbf(r)D^*\ZA\int _0^{T-r}H(T-r-s)R_-(s)\zbeta\ZD s\,\ZD \Gamma\,\ZD r=\\
\ZLA{PriCOMPO}&=\intT\int_\Gamma \zbf(r)\left \{ D^* A\left [    \ZA^{-1} R_-(T-r)\zbeta+ \int_0^{T-r}H(T-r-s)\ZA^{-1} R_-(s)\zbeta\ZD s \right ]\right \}\ZD\Gamma\,\ZD r\,.
\end{align}  
The bracket  is the  solution of~(\ref{equa:viscoELASICO}) with $\zbf=0$, $\zbF=0$ and initial condition $\zbw(0)=0$, $\zbw'(0)=\zbeta\in L^2(\ZOMq)$ so that the brace belongs to $L^2(0,T;L^2(\Gamma))$ thanks to the direct inequality.

 We have $(\zbxi,\zbeta)\perp R_V(T)$ when~(\ref{SecCOMPO}) and~(\ref{PriCOMPO}) sum to zero. 
 Taking into account that  the previous computation holds for every $f\in\mathcal{D}(\Gamma\times(0,T))$ (which is dense in $L^2(0,T;L^2(\Gamma))$   we get 
 
\[
D^*A\left \{\zbu(t)+\intt H(t-r)\zbu(r)\ZD r\right \}=0\,,\qquad \zbu(t)=R_+(t)\zbxi+\ZA^{-1} R_-(t)\zbeta\,.
\]
Let
\[
\zbpsi(t)=\zbu(t)+\intt H(t-r)\zbu(r)\ZD r\,.
\]
The function $\zbpsi(t)$ is the solution of the Volterra integral equation
\[
\zbpsi(t)=\zbu(t)+\ZA^{-1}\intt L(t-s)\zbpsi(s)\ZD s\,.
\]
We compare with~Remark~\ref{rema:perORToGOn}
 and we see that  $\zbpsi(t)$ solves
\begin{equation}\ZLA{Eq:ortog}
\zbpsi''=\zL \zbpsi +b\zbpsi(t)+\intt K(t-s)  \zbpsi(s) \ZD s \,,\qquad   \left\{\begin{array}{l}
 \zbpsi(0)= \zbxi \\
 \zbpsi'(0)= \zbeta \\  \zbpsi=0\ {\rm on}\ \partial\ZOMq\,.
 \end{array}\right.
\end{equation}
Note that the operator $\zL$ is not in the memory term and that $K(t)$ is square integrable.

In conclusion, we can state:
 \begin{Theorem}\ZLA{Theo:caraORTOG}
 We have $(\zbxi,\zbeta)\perp R_V(T)$ if and only if the solution $\zbpsi(x,t)$ of the problem~(\ref{Eq:ortog})  
satisfy the condition
\begin{equation}
\ZLA{eq:diOrtog}
D^*A \zbpsi=0\,.
\end{equation}
 \end{Theorem}
 As we noted, $D^*A=-\mathcal{T}$ when applied to the elements of ${\rm dom}\,A$, in particular when $\zbxi\in\mathcal{D}(\ZOMq)$ and $ \zbeta\in\mathcal{D}(\ZOMq)$. The direct inequality implies that $\mathcal{T}\in \mathcal{L}\left (H^1_0(\ZOMq)\times L^2(\ZOMq),L^2(0,T;L^2(\Gamma)\right )$ and so $\mathcal{ T}$ is the continuous extension of   $D^*A$ to $H^1_0(\ZOMq)\times L^2(\ZOMq)$.
 Hence  we can also state:
\begin{Corollary}\ZLA{Coro:caraORTOG}
We have $(\zbxi,\zbeta)\perp R_V(T)$ if and only if the solution $\zbpsi(x,t)$ of the problem~(\ref{Eq:ortog}) satisfies
\[
\mathcal{T}\zbpsi=0\ {\rm on}\ \Gamma
\,.
\]
\end{Corollary}

\subsubsection{Step 2-2: Fourier expansion and regularity of the elements of $\left  [R_V(T+\ZEP ) \right ]^\perp$}

We recall that $A$ is a selfadjoint operator  with compact resolvent so that its spectrum is a sequence  of eigenvalues and there exists an orthonormal basis $\{\zbphi_n\}$ whose elements are eigenvectors of $A$. We denoted $-\zl_n^2$ the eigenvalue of $\zbphi_n$ (eigenvalues can be repeated, each one a finite number of times) and we assumed  that the order has been chosen so to have $\zl_n\leq \zl_{n+1}$.
The bases of respectively $H^1_0(\ZOMq)$ and $H^{-1}(\ZOMq)$ which correspond to the orthonormal basis $\{\zbphi_n\}$ of the eigenvectors of $A$ in $L^2(\ZOMq)$ are respectively $\{\zbphi_n/\zl_n\}$ and $\{\zl_n\zbphi_n\}$. 

 Let $  (\zbxi,\zbeta)\in \left  [R_V(T+\ZEP) \right ]^\perp$,  $(\zbxi,\zbeta)\in H^1_0(\ZOMq)\times L^2(\ZOMq) $. We can expand in series of the eigenvectors and we get
\begin{equation}\ZLA{eq:expaXIedETA}
\zbxi(x)=\ZSUno \zbphi_n(x) \xi_n  \,,\quad \zbeta(x)=\ZSUno \zbphi_n(x)\eta_n\,,\qquad \{\zl_n\xi_n\}\in l^2\,,\ \{\eta_n\}\in l^2\,.
\end{equation}
First  we prove the following result (which   holds also with $\ZEP=0$):
\begin{Lemma}\ZLA{lemmaSERIEInfi}
Assume $R_E(T  )=H^{-1}(\ZOMq)\times L^2(\ZOMq) $ and   $\ZEP\geq 0$.
Let  $(\zbxi,\zbeta)\in \left  [R_V(T+\ZEP ) \right ]^\perp$ and $(\zbxi,\zbeta)\neq 0$. Then at least one of the series in~(\ref{eq:expaXIedETA})  is  not a finite sum.
\end{Lemma}
\zProof Let us denote $\zbpsi(t)$ the solution of Eq.~(\ref{Eq:ortog}) with $\zbpsi(0)=\zbxi$ and $\zbpsi'(0)=\zbeta$ and let us assume that both the expansions~(\ref{eq:expaXIedETA}) are finite sums, of $N$ terms at most. We expand also the solution $\zbpsi(t)$, 
\[
\zbpsi(t)=\sum _{n=1}^N  \zbphi_n(x)\left [\psi^0_n(t)\xi_n+\psi^1_n(t)\eta_n\right ]
\]
where both $\psi^0_n(t)$ and $\psi^1_n(t)$ solve the scalar equation
\[
\psi_n''(t)=-\zl_n^2\psi_n(t)+b\psi_n(t)+\intt K(t-s)\psi_n(s)\ZD s
\]
with initial conditions respectively
\[
\psi^0_n(0)=1\,,\ (\psi^0_n)'(0)=0\,,\quad \psi^1_n=0\,,\ ( \psi^1_n)'(0)=1\,.
\]
The orthogonality condition is
\begin{equation}\ZLA{eq:OrToCondiFini}
0=\mathcal{T}\zbpsi(t)=\sum _{n=1}^N \left (\mathcal{T} \zbphi_n(x)\right )\left [\psi^0_n(t)\xi_n+\psi^1_n(t)\eta_n\right ]
\end{equation}
This sum cannot have only one nonzero addendum.
To see this, let 
  $\xi_n=0$ and  $\eta_n=0$ for $n\neq n_0$ and either $\xi_{n_0}\neq 0$ or $\eta_{n_0}\neq 0$ (or both). Computing either~(\ref{eq:OrToCondiFini}) or its derivative with $t=0$ we get
\begin{equation}\ZLA{eq:tracciaNULLA}
 \left (\mathcal{T} \zbphi_{n_0}(x)\right )=0\,.
\end{equation}
It is a fact that the equality~(\ref{eq:tracciaNULLA}) cannot hold if $\Gamma$ has been chosen in such a way that    the elastic system~(\ref{equa:puroELASICO}) is controllable at some time $T$. 
The proof, based on the inverse inequality of Eq.~(\ref{equa:puroELASICO}), is the same as that with $d=1$ and can be found in~\cite{TAOcorrect}   see also~\cite[Lemma~4.3]{PandLIBRO}.

Even more: \emph{the nonzero elements of~(\ref{eq:OrToCondiFini}) must correspond to at least  two different eigenvalues.} In fact, Eq.~(\ref{eq:OrToCondiFini})   computed with $t=0$ gives 
\[
\mathcal{T}\left ( \sum_{n=1}^N \xi_n\zbphi_n(x)  \right )=0
\]
and if every $\zbphi_n$ correspond to the same eigenvalue $-\zl^2$ then $ \sum_{n=1}^N \xi_n\zbphi_n(x)$ is an eigenvector of $A$. As we noted, the  equality  is impossible.

Analogous argument if the derivative of both the sides of~(\ref{eq:OrToCondiFini}) is computed with $t=0$.

Now we prove that if the sum is finite then we can reduce ourselves to have a sum of only one eigenvalue, and we proved that this  is not possible. In fact, computing the second derivatives of~(\ref{eq:OrToCondiFini})  we get
  
\begin{equation}\ZLA{eq:OrToCondiFiniDeri}
\sum _{n=1}^N \zl_n^2 \left (\mathcal{T}   \zbphi_n(x)\right ) \left [\psi^0_n(t)\xi_n+\psi^1_n(t)\eta_n\right ] =0\,.
\end{equation}
In fact, the other terms appearing in the computation of the derivative sum to zero thanks to~(\ref{eq:OrToCondiFini}).
 
We subtract~(\ref{eq:OrToCondiFini}), multiplied with $\zl_N^2$,  and~(\ref{eq:OrToCondiFiniDeri}) (as in~\cite{PandMADRID})  so to obtain a new sum like~(\ref{eq:OrToCondiFini}) but with at most $N-1$ terms:
\[
\sum _{n=1}^{N-1}\left (\mathcal{T}\zbphi_n(x)\right )\left \{
\psi_n^0(t)\left [\zl_n^2-\zl_N^2\right ]\zbxi_n+
\psi_n^1(t)\left [\zl_n^2-\zl_N^2\right ]\zbeta_n
\right \}=0
\]

 After a finite number of iteration we remain with terms which correspond to the same eigenvalue (possibly, only one term)  and we have seen that this is impossible.\zdia

Now we have this information, that  
 at least one of the series~(\ref{eq:expaXIedETA}) is not   a finite sum.  
So,  the orthogonality condition~(\ref{eq:OrToCondiFini}) has to be replaced with 
\begin{equation}
\ZLA{condiORTOperSER}
0=\mathcal{T}\zbpsi(t)=\sum _{n=1}^{+\ZIN} \left (\mathcal{T} \zbphi_n(x)\right )\left [\psi^0_n(t)\xi_n+\psi^1_n(t)\eta_n\right ]
\end{equation}
(the exchange of $\mathcal{T}$ and the series is justified by the direct inequality).
   We are going to prove that also this case is impossible  and so it must be $\zbxi=0$ and $\zbeta=0$, but now we need the assumption that the purely elastic system~(\ref{equa:puroELASICO}) is controllable at time $T $ and that $(\zbxi,\zbeta)\in\left  [R_V(T+\ZEP ) \right ]^\perp$ with $\ZEP>0$; i.e. that the equality~(\ref{condiORTOperSER}) holds in $L^2(0,T+\ZEP;L^2(\Gamma))$.  The condition $\ZEP>0$ is used in the proof of the following lemma:

\begin{Lemma}\ZLA{lemma:REGOcoeff}
Let $R_E(T)=H^{-1}(\ZOMq)\times L^2(\ZOMq)$ and let $\ZEP>0$.
If $(\zbxi ,\zbeta )\in\left  [R_V(T+\ZEP) \right ]^\perp $   then
\begin{equation}\ZLA{eq:StiMExiNetaN}
 \xi_n=\frac{\tilde\xi_n}{\zl_n^3}\,,\quad \eta_n=\frac{\tilde \eta_n}{\zl_n^2}\,, \quad \{\tilde \xi_n\}\in l^2\,,\ \{\tilde\eta_n\}\in l^2\,.
\end{equation}
\end{Lemma}
 \zProof 
 We expand in series of the eigenfunctions $\zbphi_n$ the solution $\zbu(t)$ of the purely elastic problem~(\ref{equa:puroELASICO})-(\ref{condiBOrdoELASTICO})  when $\zbF=0$, $\zbu_0=0$ and $\zbu_1=0$. 
Using~(\ref{equa:puroELASICO}) and~(\ref{eq:INTEperPARTI}) (recall that $-\zbf$ was renamed $\zbf$) we get:
\[
\zbu(t)=\ZSUno\zbphi_n(x) u_n(t)\,,\qquad u_n'' =-\zl_n^2 u_n +\int_{\Gamma}  \left(\mathcal{T}\zbphi_n\right )\cdot \zbf(t)\ZD\Gamma\,.
\]
Hence,
\[
u_n(t)=\frac{1}{\zl_n}\intt \sin\zl_n(t-s)\left [\int_{\Gamma}  \left( \mathcal{T}\zbphi_n\right )\cdot \zbf(s)\ZD\Gamma\right ]\ZD s
\]
and we have the following expansions of $\zbu(t)$ and $\zbu'(t)$ (we compute with $\zbu_0=0$, $\zbu_1=0$ and $\zbF=0$):
 \begin{align*}
 \zbu(t)&=-\ZSUno\zbphi_n(x) \intt \int_\Gamma \left (\frac{\mathcal{ T}\zbphi_n}{\zl_n}\sin\zl_n(t-s)\right )\cdot \zbf(s)\ZD\Gamma\, \ZD s\,,\\
  \zbu'(t)&=-\ZSUno\left (\zl_n\zbphi_n(x)\right ) \intt \int_\Gamma
   \left (\frac{\mathcal{ T}\zbphi_n}{\zl_n}\cos\zl_n(t-s)\right ) \cdot \zbf(s)\ZD\Gamma\, \ZD s\,
 \end{align*}
 So, controllability is equivalent to the surjectivity of the  map 
\begin{align*}
\zbf\mapsto&\left \{\left [ \intT \int_\Gamma
   \left (\frac{\mathcal{ T}\zbphi_n}{\zl_n}\cos\zl_ns\right ) \cdot \zbf(T-s)\ZD\Gamma\, \ZD s\right ]\,,\right.\\
 &\left.  \left [
    \intT \int_\Gamma \left (\frac{\mathcal{ T}\zbphi_n}{\zl_n}\sin\zl_ns\right )\cdot \zbf(T-s)\ZD\Gamma\, \ZD s
   \right ]
   \right \}\in l^2\times l^2\,.
\end{align*}
Here $l^2=l^2(\mathbb{N})$, $\mathbb{N}=1,\, 2,\,\dots$.
This transformation is continuous since $\zbf\mapsto(\zbu'(T),\zbu(T))$ is continuous from $L^2(0,T;L^2(\Gamma))$ to $H^{-1}(\ZOMq)\times L^2(\ZOMq)$.  

 As usual with Fourier series, it is convenient to introduce 
\[
\zl_{-n}^2=-\zl_n^2\,,\  \zbphi_{-n}=\zbphi_n
\]
and we see that controllability of the purely elastic system is equivalent to surjectivity of the following operator  (here $\mathbb{Z}'=\mathbb{Z}\setminus\{0\}$ and $l^2=l^2(\mathbb{Z}')=l^2(\mathbb{Z}';\mathbb{C})$)
 
\begin{equation}\ZLA{eq:opeMO}
\Mo\in \mathcal{L}\left (L^2\left (0,T;L^2(\Gamma)\right ) , l^2 \right )\,:\quad \Mo \zbf=\left \{\intT \int_\Gamma
   \left (\frac{\mathcal{ T}\zbphi_n}{\zl_n}e^{i\zln s}\right )\cdot \zbf(T-s)\ZD\Gamma\, \ZD s\right \} 
\end{equation}
 (see Lemmas 4.6 and 5.1 in~\cite{PandLIBRO}): 
\emph{controllability of problem~(\ref{equa:puroELASICO}) and~(\ref{condiBOrdoELASTICO}) is equivalent to the surjectivity of the bounded operator $\Mo$ in~(\ref{eq:opeMO}) from $L^2(0,T;L^2(\Gamma))$ to $l^2(\zzzP;\mathbb{C})$.} In turn, this is equivalent to the fact that the sequence 
\[
\left \{  \left (\zbPsi_ne^{i\zl_n t}\right)\right  \}_{n\in\zzzP}\quad \mbox{where} \quad \zbPsi_n=  \frac{1}{\zl_n}\mathcal{ T} \zbphi_n 
\]
is a Riesz sequence in $L^2(0,T;L^2(\Gamma))$  i.e. it can be transformed to an orthonormal sequence using a linear bounded and boundedly invertible transformation.

We need   the following pieces of information  (see~\cite[ch.~3]{PandLIBRO} for details on the Riesz sequences):
\begin{Lemma}\ZLA{lemma:piecesINFO}
The following properties hold:
\begin{enumerate}
\item\ZLA{lemma:piecesINFO1} if $\{e_n\}$ is a Riesz sequence in a Hilbert space $H$ then $\sum \zaa_n e_n$ converges in $H$ if and only if $\{\zaa_n\}\in l^2$;
\item\ZLA{lemma:piecesINFO2}  if $\left \{k_n e^{i\zl_n t}\right \}$ is a Riesz sequence in $L^2(0,T;H)$ ($H$ is a Hilbert space and $k_n\in H$) and if \ $\sum  \zaa_n k_n e^{i\zl_n t}$ converges in $L^2(0,T+\ZEP;H)$ to an $H^1$ function then $\zaa_n=\ZDE_n/\zl_n$ and $\{\ZDE_n\}\in l^2$. \emph{This result requires  $\ZEP>0$.}
\end{enumerate}
\end{Lemma}

We continue the proof of Lemma~\ref{lemma:REGOcoeff}: we    go back to examine the orthogonality condition $\mathcal{T}\zbpsi(t)=0$ which can be written as
\begin{equation}\ZLA{eq:FormAzN}
\ZSUno \zbPsi_n Z_n(t)=0\,,\quad Z_n(t)=\zl_n\left [\psi^0_n(t)\xi_n+\psi^1_n(t)\eta_n\right ]\,.
\end{equation}
Hence $Z_n(t) $ solves
\begin{equation}
\ZLA{eqZnMac}
Z_n''=-\zl_n^2Z_n+bZ_n+\intt K(t-s)Z_n(s)\ZD s\,,\qquad Z_n(0)=\zl_n\xi_n\,,\quad Z_n'(0)=\zl_n\eta_n\,.
\end{equation}
So we have
\begin{align}
\nonumber
Z_n(t)=& (\zl_n  \xi_n) \cos\zl_n t+ \eta_n  \sin\zl_n t+\frac{b}{\zl_n} \intt\sin\zl_n(t-s)Z_n(s)\ZD s+\\
 & +\frac{1}{\zl_n}\intt \sin\zl_n(t-s)\ints K(s-r)Z_n(r)\ZD r\,\ZD s
 \ZLA{eqZnMacInteg} 
\end{align}
(note that $\xi_n$ and $\eta_n$ are real numbers and that $\{\zl_n\xi_n\}\in l^2$ because $\xi\in H^1_0(\ZOMq)$).

Gronwall inequality shows that the sequence of continuous functions $\{Z_n(t)\}$ is uniformly bounded on compact intervals.

We have
\[
(\zl_n  \xi_n) \cos\zl_n t+ \eta_n  \sin\zl_n t=e^{i\zl_nt} c_n+ e^{-i\zl_n t}\bar c_n\,,\qquad c_n =\frac{1}{2}\left [ \zl_n \xi_n-i  \eta_n \right] \,.
\]
We  introduce the notations
 
\[
U_n=e^{i\zl_nt} c_n+ e^{-i\zl_n t}\bar c_n\,,\quad S_n=\sin\zl_n t\,,\quad C_n=\cos\zl_n t\,,\quad H_n=bS_n+S_n*K\,.
\]
 
Three steps of Picard iteration applied to~(\ref{eqZnMacInteg}) give
\[
Z_n=U_n+\frac{1}{\zl_n }
 U_n*\left [
  \sum_{\nu=1}^4\frac{1}
 {\zl_n^{\nu-1}}
 H_n^{*(\nu )} 
 \right ]
 +\frac{1}{	\zl_n^5} M_n 
\]
where $M_n=M_n(t)$ are continuous functions, and the sequences  $\{M_n(t)\}$,  $\{(1/\zl_n)M'_n(t)\}$ are   bounded on bounded intervals. In fact,
\[
M_n =U_n*\left [\sum _{\nu=5}^7 \frac{1}{\zl_n^{\nu-5}}H^{*(\nu)}\right ]
+\frac{1}{\zl_n^3}H_n^{*(8)}*Z_n\,.
\]
 
The orthogonality condition~(\ref{eq:FormAzN})  takes the form
\begin{equation}\ZLA{condiOrtoEsteConVO}
-\ZSUno\zbPsi_nU_n=\ZSUno \zbPsi_n \frac{1}{\zl_n}U_n*\left [
\sum _{\nu=1}^4 \frac{1}{\zl_n^{\nu-1}}H_n^{*(\nu )} 
\right ]+\ZSUno\zbPsi_n\frac{1}{\zl_n^5}M_n  \,.
\end{equation}
In fact, we can distribute the series since every one of the obtained series converges since {\bf 1)} the last series converges because
the estimate~(\ref{eq:STImasintAUTOV}) implies that
\[
\ZSUno \frac{1}{\zl_n^4}<+\ZIN\,.
\]
{\bf 2)} the previous series converge since $S_n*U_n$ is a linear combination of exponentials $e^{i\zl_n t}$  multiplied 
with a polynomial of degree at most $1$ and our assumption on the time  $T$ implies that $\zbPsi_ne^{i\zl_nt}$ is a Riesz sequence.

The right hand side of~(\ref{condiOrtoEsteConVO}) is an $H^1$ function on \emph{every   interval.}  In fact, computing the derivative of the right hand side termwise we get the following  sum of   $L^2$-convergent series, where $V_n=U_n' =e^{i\zl_nt}\left (i\zl_n c_n\right )-
e^{-i\zl_nt}\left (i\zl_n \bar c_n\right )
$:
\begin{align} 
\nonumber&\ZSUno \zbPsi_n\frac{1}{\zl_n}\left (c_n+\bar c_n\right )\left [
 \sum_{\nu=1}^4\frac{1}
 {\zl_n^{\nu-1}}
 H_n^{*(\nu )} 
\right ]  +\\
\ZLA{eq:metaDERIvata}
& \ZSUno \zbPsi_n\left (\frac{1}{\zl_n}V_n\right )*\left [
 \sum _{\nu=1}^4 \frac{1}{\zl_n^{\nu-1} }H^{*(\nu )} 
 \right ] +\ZSUno\zbPsi_n \frac{1}{\zl_n^5}M_n'  \,.
\end{align}

Convergence of the series is seen by an argument similar to the one used above,   thanks to~(\ref{eq:STImasintAUTOV}) and to  the fact that $\{(1/\zl_n)M_n'(t)\}$ is a bounded sequence. So, we can apply Lemma~\ref{lemma:piecesINFO} and we conclude
\begin{equation}\ZLA{eq:IniDatAinterME}
c_n=\frac{1}{\zl_n}\tilde c_n\,,\  \{\tilde c_n\}\in l^2\ \  {\rm hence}\ \
\zl_n\xi_n=\frac{1}{\zl_n}\ZDE_n\,, \  \eta_n=\frac{1}{\zl_n}\ZSI_n\,,\  \{\ZDE_n\}\in l^2\,,\  \{\ZSI_n\}\in l^2\,.
\end{equation}

We equate the derivative of the left hand side of~(\ref{condiOrtoEsteConVO}) and~(\ref{eq:metaDERIvata}) and we find the equality
\begin{align}
-\nonumber& \ZSUno\zbPsi_n V_n =
\ZSUno \zbPsi_n\frac{1}{\zl_n}\left (\tilde c_n+\bar {\tilde c}_n\right )\left [
 \sum_{\nu=1}^4\frac{1}
 {\zl_n^{\nu-1}}
 H_n^{*(\nu )} 
\right ] +
\\
\ZLA{eq:normalitaALprimoPASSO}&+
 \ZSUno \zbPsi_n\left (\frac{1}{\zl_n}V_n\right )*\left [
 \sum _{\nu=1}^4 \frac{1}{\zl_n^{\nu-1} }H^{*(\nu )}
 \right ]+\ZSUno\zbPsi_n \frac{1}{\zl_n^5}M_n' \,.
\end{align}  
Using
\[
c_n+\bar c_n=\frac{1}{\zl_n}\left ( \tilde c_n+\bar{\tilde c}_n   \right )\,,\quad 
\frac{1}{\zl_n} V_n(t)=\frac{i}{\zl_n}\left [
e^{i\zl_n t}\tilde c_n- e^{-i\zl_n t}\overline{\tilde c }_n
\right ]
\]
we see differentiability of the first and second series on the right hand side. Using again the estimate~(\ref{eq:STImasintAUTOV}) we see that also  the last  series is differentiable, since
  the explicit expression of $M_n'(t)$   is 
  \begin{align*} 
\nonumber &\frac{1}{\zl_n^5} M_n'(t)=  \frac{1}{\zl_n^5}
\left [(c_n+\bar c_n)\sum _{\mu=5}^7\frac{1}{\zl_n^{\mu-5}}H_n^{*(\mu)}+
V_n'*
\left (
\sum _{\mu=5}^7\frac{1}{\zl_n^{\mu-5}}H_n^{*(\mu)}
\right ) 
\right ]+\\
&+\frac{1}{\zl_n^7}\left (C_n+C_n*K\right )*H_n^{*(7)}*Z_n\,.
  \end{align*}
  In conclusion, the left hand side of~(\ref{eq:normalitaALprimoPASSO}) is an $H^1$ function on every interval $[0,T]$.
Applying again Lemma~\ref{lemma:piecesINFO}  we get
\[
\tilde c_n=\frac{1}{\zl_n}\hat c_n\,,\ \{\hat c_n\}\in l^2\ \ \mbox{so that}\ \ 
\zl_n\xi_n=\frac{1}{\zl_n^2} \left (\hat c_n+\bar{\hat c}_n\right )\,,\quad \eta_n=\frac{-i}{\zl_n^2}\left (
\bar{\hat c}_n-\hat c_n\,,
\right )
\]
as we wanted.\zdia

Note that we proved also this result:
\begin{Theorem}\ZLA{teoDERIVABI}
The first and the second derivatives of the series in~(\ref{eq:FormAzN}) can be  computed termwise.
\end{Theorem}

\subsubsection{Step 2-3: end of the proof}
We use Theorem~\ref{teoDERIVABI}:   the second derivative of~(\ref{eq:FormAzN}) computed termwise gives the equality
 
\[
 \ZSUno \zbPsi_n Z''_n(t)=-\ZSUno \zbPsi_n \left (\zl_n^2 Z_n(t)\right )-\intt M(t-s)\left[\ZSUno \zbPsi_n \left (\zl_n^2 Z_n(s)\right )\right ]\ZD s=0
\]
 so that we have also
\begin{equation}
\ZLA{eq:ConDIortogINTERMzPRIMO}
\ZSUno \zbPsi_n \zl_n^2 Z_n(t) 
=\ZSUno\left (\mathcal{T}\zbphi_n(x)\right )\left [\psi_n^0(t)\frac{\tilde \xi_n}{\zl_n}+\psi_n^1(t)\tilde\eta_n\right ]= 
0\,.
\end{equation}

Let $N_1$ be the first index such that $\xi_{N_1}^2+\eta_{N_1}^2\neq 0$.
Combining~(\ref{eq:FormAzN}) and~(\ref{eq:ConDIortogINTERMzPRIMO}) we get the new equality
\begin{align*}
&\sum _{n=N_1}^{+\ZIN} \zbPsi_n \left (\zl^2_{N_1}-\zl_n^2\right ) Z_n(t)=
\sum _{n=N_2}^{+\ZIN} \zbPsi_n \left (\zl^2_{N_1}-\zl_n^2\right ) Z_n(t)
 \\
&=
\sum _{n=N_2}^{+\ZIN} \left (\mathcal{T}\zbphi_n(x)\right )
\left [
\left(
\frac{\zl_{N_1}^2}{\zl_n^2}-1
\right )
\frac{\tilde \xi_n}{\zl_n }\psi_n^0(t) +\left (\frac{\zl_{N_1}^2}{\zl_n^2}-1\right )
\tilde\eta_n
\psi_n^1(t) 
\right ]= 
0 
\end{align*}
and $N_2>N_1$. 
So,  the element $(\zbxi^1,\zbeta^1)\in H^1_0(\ZOMq)\times L^2(\ZOMq)$ whose Fourier coefficients are
\[
\xi_n^1=\left(
\frac{\zl_{N_1}^2}{\zl_n^2}-1
\right )
\frac{\tilde \xi_n}{\zl_n }\,,\quad \eta_n^1=\left (\frac{\zl_{N_1}^2}{\zl_n^2}-1\right )
\tilde\eta_n\,,\qquad \{\zl_n\tilde\xi_n\}\in l^2\,,\ \{\tilde\eta_n\}\in l^2
\]
is a second element in $\left  [R_V(T+\ZEP) \right ]^\perp$. 
 
 The elements $(\zbxi,\zbeta)$ and $(\zbxi^1,\zbeta^1)$ \emph{are not colinear,}  since the Fourier coefficient of index $N_1$
is nonzero for the first pair, while it is zero for the second one.

Now we use the fact that $(\zbxi_1,\zbeta_1)\in H^1_0(\ZOMq)\times L^2(\ZOMq)$ and we iterate the procedure. We remove the first nonzero Fourier coefficient of $(\zbxi_1,\zbeta_1)$ and so we get  an element $(\zbxi_2,\zbeta_2)\in H^1_0(\ZOMq)\times L^2(\ZOMq)$ which belongs to $\left  [R_V(T+\ZEP) \right ]^\perp$ and \emph{which is linearly independent from the previous ones.} If the series of $(\zbxi,\zbeta)$ is not a finite sum then the procedure can be repeated and we get the \emph{contradictory statement} that $\left  [R_V(T+\ZEP) \right ]^\perp$ has infinite codimension.

We combine this fact with Lemma~\ref{lemmaSERIEInfi} and we get $ \zbxi=0$, $\zbeta=0$ as we wanted to prove.

\end{document}